\documentclass[reqno,12pt]{amsart}
\usepackage{amscd,amssymb}
\usepackage{color}

\theoremstyle{plain}
\newtheorem{Thm}[subsection]{Theorem}
\newtheorem{Cor}[subsection]{Corollary}
\newtheorem{Lem}[subsection]{Lemma}
\newtheorem{Prop}[subsection]{Proposition}
\newtheorem{Conj}[subsection]{Conjecture}
\newtheorem{Claim}[subsection]{Claim}

\theoremstyle{definition}
\newtheorem{Def}[subsection]{Definition}

\theoremstyle{remark}

\newtheorem{Rem}[subsection]{Remark}

\newtheorem{conj}{Conjecture}

\numberwithin{equation}{section}
\renewcommand{\rm}{\normalshape}

\newif\ifShowLabels
\ShowLabelstrue
\newdimen\theight
\def\TeXref#1{%
    \leavevmode\vadjust{\setbox0=\hbox{{\tt
        \quad\quad  {\small \rm #1}}}%
    \theight=\ht0
    \advance\theight by \lineskip
    \kern -\theight \vbox to
    \theight{\rightline{\rlap{\box0}}%
    \vss}%
    }}%

\ShowLabelsfalse

\renewcommand{\sec}[2]{\section{#2}\label{S:#1}%
    \ifShowLabels \TeXref{{S:#1}} \fi}
\newcommand{\ssec}[2]{\subsection{#2}\label{SS:#1}%
    \ifShowLabels \TeXref{{SS:#1}} \fi}

\newcommand{\refs}[1]{Section ~\ref{S:#1}}

\newcommand{\reft}[1]{Theorem ~\ref{T:#1}}

\newcommand{\refc}[1]{Corollary ~\ref{C:#1}}

\newcommand{\refe}[1]{\eqref{E:#1}}
\newcommand{\refco}[1]{Conjecture ~\ref{Co:#1}}

\newenvironment{thm}[1]%
    { \begin{Thm} \label{T:#1}  \ifShowLabels \TeXref{T:#1} \fi }%
    { \end{Thm} }

\renewcommand{\th}[1]{\begin{thm}{#1} \sl }
\renewcommand{\eth}{\end{thm} }

\newenvironment{lemma}[1]%
    { \begin{Lem} \label{L:#1}  \ifShowLabels \TeXref{L:#1} \fi }%
    { \end{Lem} }
\newcommand{\lem}[1]{\begin{lemma}{#1} \sl}
\newcommand{\elem}{\end{lemma}}

\newenvironment{propos}[1]%
    { \begin{Prop} \label{P:#1}  \ifShowLabels \TeXref{P:#1} \fi }%
    { \end{Prop} }
\newcommand{\prop}[1]{\begin{propos}{#1}\sl }
\newcommand{\eprop}{\end{propos}}

\newenvironment{corol}[1]%
    { \begin{Cor} \label{C:#1}  \ifShowLabels \TeXref{C:#1} \fi }%
    { \end{Cor} }
\newcommand{\cor}[1]{\begin{corol}{#1} \sl }
\newcommand{\ecor}{\end{corol}}

\newenvironment{defeni}[1]%
    { \begin{Def} \label{D:#1}  \ifShowLabels \TeXref{D:#1} \fi }%
    { \end{Def} }
\newcommand{\defe}[1]{\begin{defeni}{#1} \sl }
\newcommand{\edefe}{\end{defeni}}

\newenvironment{remark}[1]%
    { \begin{Rem} \label{R:#1}  \ifShowLabels \TeXref{R:#1} \fi }%
    { \end{Rem} }
\newcommand{\rem}[1]{\begin{remark}{#1}}
\newcommand{\erem}{\end{remark}}

\newenvironment{conjec}[1]%
    { \begin{Conj} \label{Co:#1}  \ifShowLabels \TeXref{Co:#1} \fi }%
    { \end{Conj} }
\renewcommand{\conj}[1]{\begin{conjec}{#1} \sl }
\newcommand{\econj}{\end{conjec}}

\newcommand{\eq}[1]%
    { \ifShowLabels \TeXref{E:#1} \fi
       \begin{equation} \label{E:#1} }
\newcommand{\eeq}{ \end{equation} }

\newcommand{\prf}{ \begin{proof} }
\newcommand{\epr}{ \end{proof} }

\newcommand\alp{\alpha}

\newcommand\del{\delta}		\newcommand\Del{\Delta}

\newcommand\lam{\lambda}		\newcommand\Lam{\Lambda}
\newcommand\sig{\sigma}

\newcommand\calA{{\mathcal{A}}}
\newcommand\calB{{\mathcal{B}}}
\newcommand\calC{{\mathcal{C}}}

\newcommand\calE{{\mathcal{E}}}

\newcommand\calH{{\mathcal{H}}}

\newcommand\calJ{{\mathcal{J}}}

\newcommand\calM{{\mathcal{M}}}
\newcommand\calN{{\mathcal{N}}}
\newcommand\calO{{\mathcal{O}}}

\newcommand\calS{{\mathcal{S}}}

\newcommand\calV{{\mathcal{V}}}
\newcommand\calW{{\mathcal{W}}}

		\newcommand\bfB{{\mathbf B}}

		\newcommand\bfG{{\mathbf G}}

		\newcommand\bfM{{\mathbf M}}

		\newcommand\bfP{{\mathbf P}}

		\newcommand\bfT{{\mathbf T}}
		\newcommand\bfU{{\mathbf U}}

		\newcommand\bfX{{\mathbf X}}

\newcommand\RR{\mathbb{R}}

\newcommand\GG{\mathbb{G}}

\newcommand\CC{\mathbb{C}}

	\newcommand\grp{{\mathfrak{p}}}

	\newcommand\gru{{\mathfrak{u}}}

\newcommand\sdp{\times \hskip -0.3em {\raise 0.3ex
\hbox{$\scriptscriptstyle |$}}} 

\newcommand\Ad{\operatorname{Ad}}

\newcommand\End{\operatorname{End\,}}

\newcommand\Hom{\operatorname {Hom}}

\newcommand\id{\operatorname{id}}
\newcommand\Id{\operatorname{Id}}

\newcommand\Sym{\operatorname{Sym}}

\newcommand\tilj{{\widetilde{j}}}

\newcommand\tilPhi{{\widetilde{\Phi}}}

\newcommand\x{\times}
\newcommand\ten{\otimes}

\newcommand{\ra}{\rangle}
\newcommand{\la}{\langle}

\renewcommand{\Id}{\text{Id}}

\newcommand\cusp{\operatorname{cusp}}
\newcommand\Ass{\operatorname{Ass}}
\newcommand\sph{\operatorname{sph}}
\newcommand\Rep{\operatorname{Rep}}
\newcommand\Sat{\operatorname{Sat}}
\newcommand\St{\operatorname{St}}

\newcommand\ocalW{\overline{\calW}}
\newcommand\tcalS{\widetilde{\calS}}
\newcommand\tcalJ{\widetilde{\calJ}}
\newcommand\tsig{\widetilde{\sig}}
\begin{document}

\title{Schwartz space of parabolic basic affine space and asymptotic Hecke algebras}
\author{Alexander Braverman and David Kazhdan}
\begin{abstract}
Let $F$ be a local non-archimedian field and  $G$ be the group of $F$-points of a split connected reductive group over $F$.
In \cite{BK-J} we define an algebra $\calJ(G)$ of functions on $G$ which contains the Hecke algebra $\calH(G)$ and is contained in  the Harish-Chandra Schwartz algebra $\calC(G)$. We consider $\calJ(G)$ as an algebraic analog the algebra $\calC(G)$.

Given a parabolic subgroup $P$ of $G$ with a Levi subgroup $M$ and the unipotent radical $U_P$ we write $X_P:=G/U_P$. Let $\calS _c(X_P)$ be the space of locally constant functions on $X_P$ with compact support and $\calS _{cusp ,c}(X_P)\subset \calS _c(X_P)$ be subspace of functions whose right shifts span a cuspidal representation of $M$.

 In this paper we study the space $\calS(X_P):=\calJ(\calS _c(X_P))$ and show that

(a) For two parabolics $P,Q$ with the same Levi subrgoups the normalized intertwining operator $I_{Q,P}$ defines and isomorphism of  $(\calJ(G),\calH(M))$-bimodules
$\calS(X_P)_{cusp}$ and $\calS(X_Q)_{\cusp}$.

(b) When $P$ is a Borel subgroup the space $\calS(X_P)$ is equal to the Schwartz space studied in \cite{BK}.

\end{abstract}
\maketitle
\sec{}{Introduction and statement of the results}
\ssec{}{Notation}Let $F$ be a non-archimedian local field with
 ring of integers $\calO$; we shall fix a generator $\kappa$ of the
 maximal ideal of $\calO$. Typically, we shall denote algebraic
 varieties over $F$ by boldface letters (e.g. $\bfG,\bfX$ etc.)
 and the corresponding sets of $F$-points -- by the corresponding
 ordinary letters (i.e. $G,X$ etc.).

Let $\bfG$ be a connected split reductive group $\bfG$ over $F$ with a Borel subgroup $\bfB$,  unipotent radical $\bfU$ be   the unipotent radical of $\bfB$ and
 $\bfT :=\bfB/\bfU$ be the Cartan torus. We fix an imbedding of
 $\bfT $ into $\bfB$.

Let $\Lam$
be the lattice of cocharacters of $\bfT$
 and $\Lam ^\vee$ be the lattice of characters of $\bfT$.
We fix a Haar measure $dg$ on $G$ and we denote by $\calH(G)$ the Hecke algebra of locally constant compactly supported functions on $G$. A choice of a Haar measure defines a structure
of a locally unital associative algebra on $\calH(G)$. As well-known the category $\calM(G)$ of  smooth $G$-modules is equivalent to the category of non-degenerate $\calH(G)$-modules.

\ssec{}{Functions on $X_P$}Let $P$ be a parabolic subgroup of $G$ with a Levi subgroup $M$ and unipotent radical $U_P$.
Let $X_P=G/U_P$. This space has a natural $G\times M$ action. Therefore the space $\calS_c(X_P)$ of locally constant compactly supported functions on $X_P$ becomes a $G\times M$ module; for convenience we are going to twist the $M$ action by the square root of the absolute value of the determinant of the $M$-action on the Lie algebra $\gru_P$ of $U_P$.

It is easy to see that $X_P$ possesses unique (up to constant) $G$-invariant measure, hence we can talk about $L^2(X_P)$. It has a natural unitary action of $G\times M$ (the action of $M$ is unitary is we twist it by $|\del_P|^{1/2}$ where $\del_P:\bfM\to \GG_m$ is the determinant of the action of $\bfM$ on the Lie algebra of $\bfU_{\bfP}$).
\ssec{}{The algebra $\calJ(G)$ and the space $\calS(X_P)$}
In \cite{BK-J} we have defined certain algebra $\calJ(G)$ of functions on $G$ which contains the Hecke algebra $\calH(G)$ which can be thought of as an algebraic version of the Harish-Chandra Schwartz space $\calC(G)$ (the definition will be  recalled in \refs{J}); in particular, $\calJ(G)$ is a smooth $G\x G$ -module. It is explained in {\em loc. cit.} that $\calJ(G)$ acts on $L^2(X_P)$ for any $P$, so we can set
$$
\calS(X_P)=\calJ(G)\cdot \calS_c(X_P),
$$
where $\calS_c(X_P)$ stands for the space of locally constant functions with compact support on $X_P$. The space $\calS(X_P)$ is a smooth $G\times M$-module.

\ssec{}{Intertwining operators -- $L^2$-version} The following result is essentially Theorem 2.1 of \cite{arthur}:
\th{parabolic-l2}
Let $P$ and $Q$ be two associate parabolics, i.e. two parabolics with the same Levi subgroup $M$. Then there exists a $G\times M$-equivariant unitary isomorphism $\Phi_{P,Q}:L^2(X_P)\widetilde{\to} L^2(X_Q)$.
These isomorphisms satisfy the following properties:
\begin{enumerate}
\item
$\Phi_{P,P}=\id$
\item
For 3 parabolic subgroups $P,Q,R$ with the same Levi subgroup $M$ we have $\Phi_{Q,R}\circ\Phi_{P,Q}=\Phi_{P,R}$.
\end{enumerate}
\eth
Note that (1) and (2) together imply that $\Phi_{Q,P}\circ\Phi_{P,Q}=\id$.

\medskip
\noindent
{\bf Warning.} The operator $\Phi_{P,Q}$ is not canonical - it depends on various choices. In what follows we are going to choose some operators $\Phi_{P,Q}$ satisfying the above requirements.

\ssec{}{Another version of Schwartz space}
For a parabolic subgroup $P$ of $G$ with chosen Levi subgroup $M$ let $\Ass(P)$ denote the set of all parabolics $Q$ containing $M$ as a Levi subgroup.
We now define another version of the $\calS'(X_P)$ of the Schwartz space of functions on $X_P$ by setting
$$
\calS'(X_P)=\sum\limits_{Q\in \Ass(P)}\Phi_{Q,P}(\calS_c(X_Q)).
$$
Below is the main result of this paper:
\th{inter-schwartz}
\begin{enumerate}
\item
We have $\calS'(X_P)_{\cusp}=\calS(X_P)_{\cusp}$. Here by $\calS(X_P)_{\cusp}$ (resp. $\calS'(X_P)_{\cusp}$) we denote the $M$-cuspidal part of $\calS(X_P)_{\cusp}$ (resp. of $\calS'(X_P)_{\cusp}$).
\item
The operator $\Phi_{P,Q}$ defines an isomorphism between $\calS(X_P)_{\cusp}$ and $\calS(X_Q)_{\cusp}$.
\end{enumerate}
\eth
It is easy to see that the Schwartz algebra $\calC(G)$ acts on $L^2(X_P)$ on the right; hence $L^2(X_P)$ is also a right module over $\calJ(M)$. The following conjecture seems very plausible but we don't know how to prove it at the moment:
\conj{levi}
\begin{enumerate}
\item
The right action of $\calJ(M)$ on $L^2(X_P)$ preserves $\calS(X_P)$ (thus making $\calS(X_P)$ into a $(\calJ(G),\calJ(M))$-bimodule.
\item
We have $\calS'(X_P)\subset \calS(X_P)$.
\item
The operator $\Phi_{P,Q}$ defines an isomorphism between $\calS(X_P)$ and $\calS(X_Q)$.

\end{enumerate}
\econj
The rest of the paper is devoted to the proof of \reft{inter-schwartz}.
\ssec{}{Examples}
Let us look at two extreme cases. First, consider the case $P=G$. In this case $\Ass(P)$ consists of 1 element and therefore $\calS'(X_P)=\calH(G)$; similarly, $\calS(X_P)=\calJ(G)$. Assertions (1), (2), (3) of \reft{inter-schwartz} now become obvious and assertion
(4) is also an easy corollary of the definition of $\calJ(G)$.

Second, let us consider the case when $P$ is a Borel subgroup $B$ of $G$. Then $M$ is a maximal split torus $T$ of $G$ and we have $\calH(T)=\calJ(T)$, so assertion (1) is obvious. On the other hand, any representation of $T$ is cuspidal, so we have $\calS(X_B)_{\cusp}=\calS(X_B)$ and $\calS'(X_B)_{\cusp}=\calS'(X_B)$. Hence assertion (4) says in this case that $\calS'(X_B)=\calS(X_B)$.
This result was conjectured in \cite{BK-J} (note that $\calS'(X_B)$ is exactly the Schwartz space of \cite{BK}); assertions (2) and (3) follow from there immediately.

\ssec{}{Acknowledgments} The project has received funding from ERC under grant agreement 669655. In addition the first-named author was partially supported by NSERC.

\sec{tempered}{Tempered representations and Harish-Chandra algebra}

\defe{} For any smooth representation $(\pi ,V)$ of $G$ and smooth vectors $v\in V, v^\star \in V^\vee$ we define the function $m_{v,v^\star}$ on $G$
 $$m_{v,v^\star} (g) := v^\star(\pi (g)v), v \in V$$ and say that
$m_{v,v^\star}$ is a  {\it matrix coefficient} of the representation $\pi$.
\edefe

The algebra $\calH(G)$ can be embedded into the  Harish-Chandra Schwartz algebra $\calC(G)$. The algebra $\calC(G)$ does not act on all smooth representations of $G$ but only on tempered ones.
Let $G$ be a locally compact group. In this subsection we only consider unitary representations $(\pi , V,(,))$ of $G$.
\defe{}
\begin{enumerate}
\item We say that a representation  $(\pi , V,(,))$ is in the closure of a representation  $(r , W,(,))$ if for any vectors
$v_1,\dots  _n\in V,$  a compact $C$ in $G$ and $\epsilon >0$
there exist vectors $w_1,\dots w_n\in W$ such that

$|(\pi (g)v_i,v_j)-(r(g)v_i,v_j)|\leq \epsilon$ for all$i, j ,1\leq i,j\leq n$ and $g\in C$.
\item
A representation $(\pi , V )$ of $G$ is {\it tempered} if it is in the closure of the regular
representation of $G$.
\item
We denote by $\calM_t(G)$ the category of tempered representations.
\item
We denote by $\hat G_t$  the set of tempered irreducible representations. As follows Schur's lemma we can consider
$\hat G_t$  as a subset of the set $\hat G$ of smooth  irreducible representations of $G$.
\end{enumerate}
\edefe
The following results are proven in \cite{CHH} and \cite{SZ1}.

\begin{Claim}{}\label{te}
\begin{enumerate}
\item

A unitary irreducible representation $(\pi , V, (, ))$ of a reductive group $G$ over a local field is tempered iff matrix coefficients of $\pi$ belong to  $L^{2+\epsilon}(G/Z(G))$ for any $\epsilon > 0.$
\item
 For any tempred representation $V$ of $G$ the action of the algebra $\calH(G)$ on $V$ extends to a continuous action of the Harish-Chandra algebra $\calC(G)$.
\item

Let $P$ be a parabolic subgroup of $G$ with a Levi group $M$ and $\sigma$  be a tempered irreducible representation of $M$. Then the unitary induced representation
$\pi_\sigma=ind^G_P(\sigma)$ is tempered.
\item
 For a generic unitary character $\chi :M\to S^1$
 the representation $\pi _{\sigma \otimes \chi}$ (which is tempered) is irreducible.
\item Any representation which is  a Hilbert space integral of tempered representations is tempered.

\end{enumerate}
\end{Claim}
Let us now discuss the Harish-Chandra algebra $\calC(G)$.
\defe{}
\begin{enumerate}
\item
Recall that for any $g\in G$, there exists a unique dominant coweight $\lam(g)$ of $\bfT$ such that $g\in \bfG(\calO)\lam(g)(\kappa )\bfG(\calO)$. We define a function $\Del$ on $G$ by
 $\Del(g):=q^{\la\lam,\rho\ra}$.

\item We say that a function $f:G\to \CC$ is a {\it Schwartz function} if

\medskip
(a) There exists an open compact subgroup $K$ of $G$ such that $f$ is two-sided  $K$-invariant.

(b) For any polynomial function $p:G\to F$ and $n>0$,  there exists a constant $C=C_{p,n}\in \RR_{>0}$ such that
$$
\Del(g)|f(g)|\leq C\ln ^{-n}(1+| p(g)|)
$$
for all $g\in G$.

\medskip
\noindent

We denote by $\calC(G)$ the space of  Schwartz functions.

\end{enumerate}
\edefe
Obviously we have an inclusion
$\calH(G)\hookrightarrow \calC(G)$ and $\calH(G)\subset \calC(G)$ is dense in the natutal topology of $\calC(G)$.

The following statements are well known    (see for example  \cite{SZ1}).
\begin{Claim}\label{cg}
\begin{enumerate}
\item $\calC (G)$ has  an algebra structure  with respect to convolution.
\item Any tempered representation of $\calH(G)$ extends to a continuous representation of $\calC(G)$.

\end{enumerate}
\end{Claim}

\begin{Lem}
\label{Temp} The natural unitary representation of $G\times M$ on the smooth part of  $L^2(X_P)$ is tempered.
\end{Lem}
\prf Since the right action of $M$ on $X_P$ is free we can write
the space $L^2(X_P)$ as a Hilbert space integral
$$
L^2(X_P)=\int \limits_{\rho \in \hat M_t}\pi (\rho )\otimes \rho ~ d\mu _M
$$
where $\mu _M$ is the Planchel measure on $\hat M_t$ and $\pi (\rho)$ is a unitary representation of $G$. It is easy to see that $\pi (\rho)=i_{GP}(\rho)$. Now Lemma \ref{Temp} follows from  Claim \ref{cg}.
\epr
\begin{Cor}{}\label{act} The natural representation of $\calH(G)$ on
the smooth part of $L^2(X_P)$ extends to a representation of $\calC(G)$.
\end{Cor}

\sec{J}{Paley-Wiener theorems and the definition of the algebra $\calJ(G)$}
\ssec{}{The Paley-Wiener theorem for $\calH(G)$}
Let $P$  be a parabolic subgroup  of $G$ with a Levi group $M$. The set  $\Psi (M)$ of unramified characters of $M$ is
equal to $\Lam ^\vee _M\otimes \CC ^\times$ where
 $\Lam ^\vee _M\subset \Lam^{\vee}$ is the subgroup of characters of $\bfT$ trivial on $\bfT \cap [\bfM,\bfM]$. So
 $\Psi (M)$ has a structure of a complex algebraic variety; the algebra of polynomial functions on $X_M$ is equal to $\CC[\Lam_M]$ where $\Lam_M$ is the lattice dual to $\Lam_M^{\vee}$. We denote by $\Psi _t(M)\subset \Psi (M)$ the subset of unitary characters.

 For any $(\sig ,V)\in\calM(M)$ we denote by $i_{GP}(\sig)$ the corresponding unitarily
 induced object of $\calM(G)$. As a representation of $G(\calO)$ this representation is
 equal to
 $\text{ind} ^{\calO}_{P(\calO)}(\sig)$.
 So for any  unramified character
$\chi:M\to \CC^*$ the space $V_\chi$ of the representation
$i_{GP}(\sig\ten \chi)$  is isomorphic to the space $V_\sig$  of the representation $i_{GP}(\sig )$ and  is independent on a choice of $\chi$.
Since  $X_M$ has a structure of an algebraic variety over $\CC$ it make sense to say that a family $\eta_\chi \in \text{End} (V_\chi ), \chi \in \Psi (M)$ is  regular or  smooth.

We denote by $Forg :\calM (G)\to Vect$ the forgetful functor,
 by $\widetilde{\calE(G)}=\{ e(\pi )\}$ the ring of endomorphisms of $Forg$ and define  $\calE(G)\subset \widetilde{\calE(G)}$ as the subring of endomorphisms $\eta_{\pi}$ such that

1) For any
 Levi subgroup $M$ of $G$ and $\sig \in Ob( \calM (M))$,
 the endomorphisms $\eta_{i_{GP}(\sig\ten \chi)}$  are regular functions of $\chi$.

2) There exists an open compact subgroup $K$ of $G$ such that $\eta_{\pi}$ is $K\times K$-invariant for every $\pi$.

By definition, we have a homomorphism
$$
PW :\calH (G)\to \calE(G),\quad f\mapsto \pi (f).
$$
The following is usually called "the matrix Paley-Wiener theorem" (cf. \cite{Ber}, Theorem 25):
\th{pw}The map $PW$ is an isomorphism.
\eth

\ssec{}{The  Paley-Wiener $\calC(G)$}

As follows from the  Claim \ref{te} the representations $\pi _{\sigma \otimes \chi}$ of $G$  belong to $\calM _t(G)$.
 for any tempered representation σ of $M$ and a unitary character $\chi$ of $M$.

Let $\calE _t(G)$ be the subring of endomorphisms $\{\eta \}$ of the forgetful functor $Forgt : \calM _t(G) \to Vect$  such that
\begin{enumerate}\item
The function $\chi \to \eta (\pi _{\sigma \otimes \chi})$
 is a smooth function of $\chi \in \Psi _t(M)$ for any
Levi subgroup $M$ of $G$ and $\sigma \in \hat M_t$.
\item
There exists an some open compact subgroup $K$ of $G$ such that
 $\eta$ is $K \times  K$-invariant.

\end{enumerate}
The following version of the matrix Paley-Wiener theorem is contained in the last section of \cite{SZ2}.

\begin{Claim}\label{tPW} The map $f  \to  \pi (f)$ defines an isomorphism between algebras $\calC(G)$ and $\calE _t(G)$.
\end{Claim}
\ssec{}{The definition of the algebra $\calJ(G)$}

Let $P$ be a parabolic subgroup with Levi group $M$. We say that an unramified character $\chi:M\to \CC^*$ is (non-strictly) positive if for any coroot $\alp$ of $\bfG$, such that the corresponding root subgroup lies in the unipotent radical $\bfU_{\bfP}$ of $\bfP$ (which in particular defines a homomorphism $\alp:F^*\to Z(M)$), we have $|\chi(\alp(x))|\geq 1$ for $|x|\geq 1$.

Let $\calE_{\calJ}(G)$ be ring of collections $\{ \eta_{\pi}\in \End_{\CC}(V)|\ \text{for tempered irreducible $(\pi,V)$}\}$ which extend to a rational function $E_{i_{GP}(\sig\ten \chi)}\in \End_{\CC}(\sig\ten \chi)$ for every tempered irreducible representation $\sig$ of $M$  and which are

a) regular on the set  of  characters $\chi$
such that $\chi^{-1}$ is (non-strictly) positive.

b) $K$-invariant for some open compact subgroup $K$ of $G$.

\noindent
As follows from the definition, we have an embedding $\calE_{\calJ}^I(G)\to \calE^I_t(G)$.

\defe{}
We define $\calJ(G)$ to be the preimage of $\calE_{\calJ}(G)$ in $\calC(G)$. Note that we have natural embeddings $\calH(G)\subset \calJ(G)\subset \calC(G)$.
\edefe
Next, let us explain certain direct sum decompositions of algebras $\calC(G)$ and $\calJ(G)$.
\defe{}
\begin{enumerate}
\item
 Let $(M ,\sig ),  (M' ,\sig ')$ be a pair of square-integrable representations of Levi subgroups of $G$. We write
  $(M ,\sig )\sim (M' ,\sig ')$ if there exists an element $g\in G$ and a unramified character $\chi \in \Psi (M)$ such that
$M'=M^g$ and $(\sig ')^g$ is equivalent to $\sig \otimes \chi$.
\item We denote by $R$ the set of equivalent classes of such  representations   $(M ,\sig )$.
\item For any     $r=(M ,\sig )$ we denote by $\calM _r(G)\subset \calM_t(G)$ the subcategory of representations in the closure of
$\bigcup\limits_{\chi \in \Psi (M)}i_{GP}(\sig \ten \chi)$ where
$P=MU_P$ is a parabolic subgoup and by $\calC_r(G)\subset \calC(G)$ the corresponding subalgebra.
\item For any $r\in R$ we define $\calJ_r(G):=\calC_r(G)\cap \calJ(G)$
\end{enumerate}
\edefe
The following statement is contained in \cite{SZ2}.
\begin{Claim}
\begin{enumerate}
\item
The subcategory $\calM _r(G)$ does depends neither on a representative $(M ,\sig )$ of $r$ nor on the choice of a parabolic $P$.
\item
$\calM _t(G)=\bigoplus\limits_{r\in R}\calM _r(G)$
\end{enumerate}
\end{Claim}
\begin{Cor}$\calJ(G)=\bigoplus\limits_{r\in R}\calJ_r(G)$
\end{Cor}
\sec{}{Intertwining operators} Until the and of this section we fix $r\in R$ and choose a representative  $(M,\sig )$  of $r$ and  a  parabolic subgroup  $P=MU_P\subset G$. We write $X_P:=G/U_P$ and denote by $dx$ a $G$-nvariant measure on $X_P$.

 For any character
 $\chi \in \Psi (M)$ we define
$(\pi _P(\chi ),V_{P,\chi}:=ind^G_P\sig _\chi$.
The restriction of the
 representation $\pi _P(\chi )$ to the subgroup $G(\calO) $  does
 not depend on $\chi$ (and is equal to to
 $ind ^{G(\calO)}_{P\cap G(\calO)}\sig _{M(\calO}$). We denote this space by $V_P$. Of course the image of $V_{P,\chi}$ in
the space of $M$-valued functions on $X_P$ depends on $\chi$. For any $f\in \calV_P,\chi \in \Psi (M)$ we denote by $f_\chi$ the corresponding function of $X_P$.
\begin{Def}
\begin{enumerate}
\item We denote $\CC (X_P)$ the space of smooth complex-valued functions on $X_P$ and $\calS _c(X_P)\subset \calS (X_P)$ be the subspace of compactly supported functions.
\item For any pair $P=MU_P,Q=MU_Q$ of associated parabolic subgroups we denote by $I_{P,Q}:\calS _C(X_P)\to \CC (X_Q)$
the {\it geometric (or non-normalized) intertwining operator} given by
$$I_{P,Q}(f)(g)=\int _{U_Q}f(gu)du$$
\end{enumerate}
\end{Def}

The following statement is contained in
Section $2$ of \cite{arthur} (Theorem 2.7); it is a slightly stronger version of \reft{parabolic-l2} .
\begin{Claim}{}\label{property}

\begin{enumerate}
\item There exists an non-empty open subset $\Psi _+(M)\subset \Psi (M)$ such that the integral $I_{P,Q}(f)=\int _{U_Q}f(gu)du$ is absolutely convergent for all $f\in \calV_\chi$.
\item For any $f\in \calV_P$ the function $I_{P,Q} (f_\chi )\in \calV_Q$
is a rational function of $\chi$.
\item The map  $I_{P,Q} :\calV_{P,\chi}\to \calV_{Q,\chi}$ is
$G$- covariant.
\item There exist  rational $\CC$-valued functions $r_{P,Q}(\chi )$ such that the operators $\Phi_{P,Q}:=r_{P,Q}I_{P,Q}$ satisfy the following:
\begin{enumerate}

\item $\Phi_{P,Q}$ is regular on the set of non-strongly elements with respect to $P$.
\item $\Phi_{P,Q}\circ \Phi_{Q,P}=Id$
\item $\Phi_{P,Q}(\chi)$ is unitary for unitary characters $\chi$.
\item For any three associate parabolic subgroups $P_i=MU_i, 1\leq i\leq 3$ we have
$\Phi_{P_1,P_2}\circ \Phi_{P_2,P_3}=\Phi_{P_1,P_3}$.
\end{enumerate}
\end{enumerate}
\end{Claim}

\ssec{}{Functions on $X_P$} We define two subspace of $L^2(X_P)$
\begin{Def}
\begin{enumerate}
\item
$
\calS(X_P):=\calJ(G)\cdot \calS_c(X_P)$
\item $
\calS '(X_P):=\sum _Q \Phi_{P,Q}(\calS_c(X_Q))$

wher the sum is over the set of parabolic subgroup $Q=MU_Q$ of $G$ associated with $P$.
\end{enumerate}
\end{Def}

\begin{Rem} Operators $\Phi_{P,Q}$ are not canonical but it is easy to see that any two choices differ by the multiplication by an regular invertible function on $\Psi (M)$ and therefore the space $  \calS'(X_P)$ is well defined.
\end{Rem}

\prop{twodef}
 There exists a surjective (but not necessarily injective) morphism of $G\times M$- modules $\calJ(G)_{U_P}\to \calS(X_P)$.

\eprop
\prf
We define a map $ \alpha _P:\calC(G)\to L^2(X_P)$
by $$ \alpha (f)(g):=\int _{U_P}f(gu)du$$
where the absolute convergence of the integral follows from Theorem 4.4.3 in \cite{SHC}. It is clear that $ \alpha _P$ factorizes through the map  $\zeta _P:\calC(G)_{U_P}\to L^2(X_P)$. It is clearly sufficient to prove the following statement.

\begin{Lem}
\begin{enumerate}
\item $\alpha _P(\calJ(G))\subset \calS(X_P)$
\item  $\alpha _P$ defines a surjection from
$\calJ(P)_{U_P}$ onto $\calS(X_P)$
\end{enumerate}
\end{Lem}
\prf Since the map $ \alpha$ commutes with the action of the algebra $\calJ(G)$ and  $\calJ(G)\calH(G)=\calJ(G)$ it is sufficient to see that  $\alpha _P(\calH(G))= \calS_c(X_P)$.
But the last claim is obviously true.
\epr
\epr

Let us explain why $\zeta$ is not necessarilly injective. Let $G=SL(2,F)$ and  $P$ be a Borel subgroup of $G$; thus $X_P$ can be naturally identified with $F^2\backslash\{ 0\}$. Let $\St$ denote the Steinberg representation of $G$. Then $\calJ(G)$  contains a direct summand isomorphic to $\End_f(\St)$ where $\End_f$ stands for endomorphisms of finite rank. It is easy to see that any homomorphism of $G$-modules from $\St$ to $L^2(F^2\backslash \{0\})$ is equal to $0$. Hence the above subalgebra must act by $0$ on $L^2(F^2\backslash \{0\})$. On the other hand, $\dim\St_{U_P}=1$ hence $\St\otimes \St_{U_P}\simeq \St$ is a non-zero subspace of $\calJ(G)_{U_P}$ which lies in the kernel of $\zeta$.

We now pass to the proof of \reft{inter-schwartz}. First we are going to discuss the explicit form of normalized intertwining operators for maximal parabolics.


\sec{}{Intertwining operators in the cuspidal  corank $1$ case}
In this section we assume that
$M\subset G$ be a Levi subgroup of semi-simple corank $1$. Then there exist two parabolic subgroups $P=MU_P$ and  $P_-=MU_-$ containing $M$.
We would like to give an explicit description of the normalized intertwining operator $\Phi_{P,P_-}$.

We denote by $M_+\subset M$ the subset of elements $m$ such that the map $u\to mum^{-1}$ contracts $U_P$ to $e$. To simplify notations we assume that $G$ is semisimple.
We can identify $\CC ^\star $ with
$\Psi (M), z\to \chi _z$
in such a way that that that $|\chi _z(m)|<1$ for
$|z|\leq 1, m\in M_+$.

Let $\sig$ be an irreducible unitary cuspidal reprentation of $M$.
We write $(\pi _z,\calV_z)$ for the representation of $G$ on $i_{GP}(\sig\ten \chi_z)$.
The following statement is contained in \cite{Ber}.

We denote by $W_\sig \subset N_G(M)/M$ the subgroup of elements
$x$ such that $\sig ^x=\sig\ten\chi$ for some $\chi\in \Psi(M)$. Since $|N_G(M)/M|\leq 2$ we see that either
$W_\sig =\{ e\}$ or $W_\sig =\{ S_2\}$.
\begin{Claim}{}\label{Ber}
\begin{enumerate}
\item If $W_\sig =\{ e\}$ then representations $\calV_z$ are irreducible for all $z\in \CC ^\star$.
\item
If $|z|\neq 1$ and the representation $(\pi _z, \calV_z)$ is reducible then $\calV_z$ has unique irreducible $G$-submodule  $W_z\subset \calV_z$ and the quotient representatikn  $\ocalW_z:=\calV_z/W_z$ is irreducible.
\item
If $|z|=1$ then either $\pi _z$ is irreducible or is the direct sum of two tempered irreducible representations $\calV_z=\calV_z^+\oplus \calV_z^-$.
\end{enumerate}
\end{Claim}
We consider now the case when $W_\sig =\{ S_2\}$. We shall write
 $I(z): \calV_{P,\chi }\to \calV_{Q,\chi}$ instead of
 $I_{P,\chi _z}$.
The following statement is proven in \cite{S}.
\begin{Claim}{}\label{S1}

 Either the geometric intertwining operator $I(z)$ is regular and invertible for all $z\in \CC ^\star$ of there exist $z_0\in \CC ^\star ,|z_0|=1$ such that $I$ has a first order pole at $z_0$.

\end{Claim}
In the first case we have $\Phi=I$.

We consider now the case when $I$ has a pole. Clearly, we can assume without loss of generality that $z_0=1$.
The next  statement is also contained in \cite{S}.
\begin{Claim}{}\label{S2}
\begin{enumerate}

\item $I(z)$ has a first order pole at $z=1$ and $(z-1)I(z)(1)=a\cdot \Id$, where $a\in \CC ^\x$.

\item There exists $c>1$ such that operators $I^{\pm 1}(z)$ are regular and invertible for
$z\notin \{ 1, c^{\pm 1}\}$.
\item
$I$ is regular at $c$, $\ker(I(c))=\calW_c$ and $I(c)$ defines an isomorphism $\ocalW_{c}\to \calW_{-,c}$.

\item $I^{-1}$  is regular at $c^{-1}, \ker(I(c^{-1}))=\calW_{c^{-1}}$ and $I(c^{-1})$ defines an isomorphism $\ocalW_{-,c^{-1}}\to \calW_{c^{-1}}$.
\end{enumerate}
\end{Claim}
Now we define $\Phi(z)=I(z)\frac{(1-z)}{(1-cz^{-1})}$.

It is now clear that $\Phi$ can be considered as a normalized intertwining operator $\Phi_{P,P_-}$; if we define a similar operator $\Phi_{P_-,P}$ then these operators  satisfy the conditions of Claim \ref{property}.

To formulate the next statement we need to slightly change the point of view.
Namely, we would like to modify $\Phi$ so that it becomes an operators from $\calV_{P,z}$ to $\calV_{P,z^{-1}}$.
For this let us choose an element $n\in N_G(M)$ so that $n^2$ belongs to the center of $M$ and $nPn^{-1}=P_-$.
Multiplying $\sig$ by some element of $\Psi(M)$ we can assume that $\Ad(n)(\sig)\simeq\sig$. Then
the right multiplication by $n$ defines an isomorphism between $\calS_c(X_P)_{\sig}$ and $\calS_c(X_{P_-})_{\sig}$
which commutes with $G$ and commutes with $M$ up to the action of $\Ad(n)$. Hence, for every $z\in \CC^*$ it defines
an isomorphism between $\calV_{P,z^{-1}}$ and $\calV_{P_-,z}$. Composing $\Phi$ with the inverse of this isomorphism
we get a (rational) isomorphism $\tilPhi(z)$ between $\calV_{P,z}$ and $\calV_{P,z^{-1}}$.

\cor{1}
\begin{enumerate}

\item $\calS '$ is the space of regular functions $f:\CC ^\star -\{ c\}\to V$
 such that the function $(z-c)f(z)$ is regular at $z=c$, $((z-c)f(z))(c)\in \calW_c$ and such that $f(z)$ is $K$-invariant for some open compact subgroup $K$ of $G$.
\item

$\calJ _\sig$ is the space of regular functions $h:\CC ^\star -\{ c\}\to \End(V)$
 such that
\begin{enumerate}

\item
the function $(z-c)h(z)$ is regular at $z=c$
\item
 $((z-c)h(z))(c)\in \Hom(\ocalW_c, \calW _c)$
\item $f(z)$ is two-sided $K$-invariant for some open compact subgroup $K$ of $G$
\item $h(z^{-1})=\tilPhi(z)h(z)\tilPhi^{-1}(z))$.
\end{enumerate}
\end{enumerate}
\ecor
\prf
The first assertion follows immediately from the definition of $\calS'$ and the construction of $\Phi(z)$.
Let us prove the 2nd assertion. First of all, it is clear that any $h$ satisfying (a)-(d) belongs to $\calJ_{\sig}$.
On the other hand, let $h$ be any element of $\calJ_{\sig}$. Then by definition it defines a rational function $h:\CC ^\star -\{ c\}\to \End(V)$
which satisfies (c) and (d) and which does not have poles when $|z|\leq 1$. Since $\tilPhi(z)$ is an isomorphism for $z\neq c^{\pm 1}, 1$ it follows that $h$ can have a pole only at $c$. Now (d) and Claim \ref{S2} imply conditions (a) and (b).
\epr
\lem{}
$\calS '=\calS$
\elem
\prf It is clear from \refc{1} that for any $j\in \calJ _\sig$ and $f\in \calS _c(X_P)_\sig $
we have $j(f)\in \calS '$. It also clear that the quotient $\calS '/\calS _c(X_P)_\sig $ is isomorphic to $\calW _c$ as a representation of $\calH$.
So to prove the equality $\calS '=\calS$ it is sufficient to find $j\in \calJ_\sig$ which is not regular at $c$. To find such $j$ choose a function $r$ in \refc{1} which satisfies the first two conditions of Corollary, has a pole at $c$ and has a second order zero at $c^{-1}$.
Now take $h(z)=r(z)+\Phi(z)r(z)\Phi^{-1}(z)$.

\epr

\sec{}{End of the proof of \reft{inter-schwartz}}

\ssec{}{}
In this Section we would like to prove \reft{inter-schwartz}.
In fact, we claim that it is enough to prove the following:

\th{emb}
 $\calS '(X_P)_{\cusp}\subset \calS (X_P)_{\cusp}$.
\eth

\ssec{}{\reft{emb} implies \reft{inter-schwartz}}
We claim that \reft{emb} implies \reft{inter-schwartz}.
First, let us prove the 2nd assertion of \reft{inter-schwartz}. For this it is enough to show that for every $\phi\in \calS(X_Q)$ we have
$\Phi_{Q,P}(\phi)\in \calS(X_P)$. By definition of $\calS(X_Q)$ we have
$$
\phi=j\cdot f
$$
where $f\in \calS_c(X_Q)$ and $j\in \calJ(G)$. Then we have
$$
\Phi_{Q,P}(\phi)=j\cdot \Phi_{Q,P}(f)\in \calJ(G)\cdot \calS'(X_P)\subset \calJ(G)\cdot \calS(X_P)=\calS(X_P).
$$

So, to complete the proof it is enough to show
that \reft{emb} implies that $\calS '(X_P)_{\cusp}=\calS (X_P)_{\cusp}$. Let us as before choose a unitary cuspidal representation $\sig$ of $M$
and let $\calS(X_P)_{\sig}$, $\calS'(X_P)_{\sig}$ be the corresponding direct summands of $\calS(X_P)_{\cusp}$, $\calS'(X_P)_{\cusp}$.
These are modules over $\CC[\Psi(M)]$ and it is clear that if we take $K$-invariant vectors for some open compact subgroup $K$ of $G$, these modules become finitely generated. Hence it is enough to show that the embedding $\calS '(X_P)_{\cusp}\subset \calS (X_P)_{\cusp}$ is surjective
in the formal  neighbourhood of every $\chi\in \Psi(M)$.

Let us first assume that $\chi$ is non-strictly negative with respect to $P$. Then by definition we have $\widehat{\calS(X_P)_{\sig,\chi}}=
\widehat{\calS_c(X_P)_{\sig,\chi}}$ where $\widehat{\calS(X_P)_{\sig,\chi}}$ and $\widehat{\calS_c(X_P)_{\sig,\chi}}$ denote the formal completions of the corresponding spaces at $\chi$. Since $\calS_c(X_P)\subset \calS'(X_P)$ the desired surjectivity follows.
Let now $\chi$ be arbitrary. Then there exists an associate parabolic $Q$ such that $\chi$ is non-strictly negative with respect to $Q$.
Since we have $\calS'(X_P)=\calS'(X_Q)$ and $\calS(X_P)=\calS(X_Q)$, it follows from the above argument that the map $\widehat{\calS'(X_P)_{\sig,\chi}}\to
\widehat{\calS(X_P)_{\sig,\chi}}$ is surjective for every $\chi$.

So, it remains to prove \reft{emb}.
\prf
We have to show that
$\Phi_{P,Q}(\calS _{c,\sig}(X_P)))\subset \calS (X_P)_{\sig}$
for all pairs  $P,Q$ of associated parabolic subgroup.

We begin by the following:
\prop{}
Let $K$ be an open compact subgroup of $G$. Let also $\calA=\CC[\Psi(M)]$. Then $\calS(X_P)^K_{\sig}$ is a reflexive $\calA$-module.
\eprop
\prf
For a finitely generated $\calA$-module $N$ let $N^{\vee}=\Hom_A(N,\calA)$. Let us compute
$(\calS(X_P)_{\sig}^K)^{\vee}$. For this let us denote by $\tcalJ(G)$ the algebra whose definition is similar to that of $\calJ(G)$ but with "non-strictly negative" replaced by "non-strictly positive". It is clear that the map $g\mapsto g^{-1}$ defines an anti-isomorphism between
$\calJ(G)$ and $\tcalJ(G)$ We denote this anti-isomorphism by $j\mapsto \tilj$. Let
$$
\tcalS(X_P)=\{ \phi\in \calS_c(X_P)|\ \tilj(\phi)\in \calS_c(X_P)\ \text{for any $\tilj\in \tcalJ(G)$}\}
$$
Similarly, we can define $\tcalS(X_P)_{\sig}$. We claim that for any $K$ as above we have
\eq{dual}
(\calS(X_P)_{\sig}^K)^{\vee}\simeq \tcalS(X_P)_{\tsig}^K\ \text{and}\ (\tcalS(X_P)_{\sig}^K)^{\vee}\simeq \calS(X_P)_{\tsig}^K,
\end{equation}
where $\tsig$ denotes the (smooth) dual representation of $\sig$.
Let us show the first isomorphism. It is clear that $(\calS_c(X_P)_{\sig}^K)^{\vee}\simeq \calS_c(X_P)_{\tsig}^K$. Hence any $\lam\in (\calS(X_P)_{\sig}^K)^{\vee}$ can be thought of as an element of $\calS_c(X_P)_{\tsig}^K$. We claim that $\tilj(\lam)\in \calS_c(X_P)_{\tsig}$.
Indeed for any $\phi\in \calS_c(X_P)_{\sig}$ we have
$$
\tilj(\lam)(\phi)=\lam(j(\phi))\in \calA
$$
since $j(\phi)\in\calS(X_P)_{\sig}$. Hence $\tilj(\lam)$ is a well-defined map from $\calS_c(X_P)_{\sig}$ to $\calA$. Conversely, assume that for some
$\lam\in\calS_c(X_P)_{\tsig}$ we have
$\tilj(\lam)\in\calS_c(X_P)_{\tsig}$. Then the same argument shows that $\lam$ is a well-defined map from $\calS(X_P)_{\sig}$ to $\calA$, i.e. an element of $(\calS(X_P)_{\sig})^{\vee}$.

Let us show the 2nd isomorphism in \refe{dual}. It is clear that we have a well-defined injective map from $\calS(X_P)_{\sig}$ to $(\tcalS(X_P)_{\tsig})^{\vee}$.

\epr
\defe{} We say that parabolic subgroups $P=MU_P$ and $Q=MU_Q$ are {\it neighbors} if there exists a parabolic subgroup $R=LU_L$
such that $U_L\subset U_P\cap U_Q$ and $M$ is subgroup of $L$ of corank $1$.
\edefe

As is well-known for any pair $P=MU_P$ and $Q=MU_Q$ there exists a chain $P_0=P, P_1,\dots ,P_l=Q$ such that $P_i, P_{i+1}$
are  neighbors for all $i,0\leq i\leq l-1$. As follows from Claim \ref{property} we have
$$
\Phi_{P,Q}=\Phi_{P,P_1}\circ \dots \circ \Phi_{P_{l-1},Q}
$$
Therefore it is clearly sufficient to show that
$\Phi_{P,Q}(\calS _c(X_P))_{\cusp})\subset \calS (X_P)_{\cusp}$ for   neighboring pairs $P,Q$.

Given such a pair, $P_L,Q_L$ be the corresponding parabolic subgroups of $L$. Let $\calS'_{P,Q,\sig}:=\Phi_{P,Q}(\calS _c(X_P)_{\sig})+\calS _c(X_Q)_{\sig}$.
Clearly, it is enough to show that $\calS'_{P,Q,\sig}\subset \calS(X_P)_{\sig}$ for any neighbouring pair $(P,Q)$.

\epr

\sec{}{Some further questions}

\ssec{}{The $\calJ$-version of the Jacquet functor}
Let us assume \refco{levi}.
Then we can define a functor
$$
r^{\calJ}_{GP}: \text{Right $\calJ(G)$-modules}\ \to \ \text{Right $\calJ(M)$-modules}
$$
by setting
\eq{}
r^{\calJ}_{GP}(\pi)=\pi\underset{\calJ(G)}\otimes \calS(X_P).
\end{equation}
It would be interesting to investigate exactness properties of this functor and compute it in some examples. Note that manifestly we have
$r^{\calJ}_{GP}(\calJ(G))=\calS(X_P)$. Also note that at least non-canonically the functor $r^{\calJ}_{GP}$ depends only on $M$ and not on $P$ (the non-canonicity comes from the fact that the operators $\Phi_{P,Q}$ are not canonically defined).

\ssec{}{The spherical part}
Let us define the spherical part of $\calS(X_P)$ by setting
\eq{}
\calS_{\sph}(X_P)=\calS(X_P)^{\bfG(\calO)\times \bfM(\calO)}.
\end{equation}
We would like to describe this space explicitly.
For this note that set-theoretically we have
$$
\bfG(\calO)\backslash X_P/\bfM(\calO)=\bfM(\calO)\backslash M/\bfM(\calO).
$$
Hence elements of $\calS_{\sph}(X_P)$ can be thought of as $\bfM(\calO)\times\bfM(\calO)$-invariant functions on $M$. Recall that the Satake isomorphism (for $M$) says that the spherical Hecke algebra $\calH_{\sph}(M)$ consisting of compactly supported $\bfM(\calO)\times\bfM(\calO)$-invariant functions on $M$ is isormophic to the complexified Grothendieck ring of the category of finite-dimensional representations of the Langlands dual group $M^{\vee}$ (considered as a group over $\CC$). We shall denote the corresponding
map $K_0(\Rep(M^{\vee})\to \calH_{\sph}(M)$ by $\Sat_M$.

Let $G^{\vee}$ denote the Langlands dual group of $G$ and let $P^{\vee}$ be the corresponding parabolic subgroup of $G^{\vee}$ with unipotent radical $U_{P^{\vee}}$. Let $\gru_{\grp^{\vee}}$ denote the Lie algebra of $U_{P^{\vee}}$; it has a natural action of $M^{\vee}$.

Now let
\eq{}
f_P=\sum\limits_{i=0}^{\infty} \Sat_M([\Sym^i(\gru_{\grp^{\vee}})]).
\end{equation}
As was discussed above, we can regard $f_P$ as a $\bfG(\calO)\times \bfM(\calO)$-invariant function on $X_P$.
\conj{spherical}
$\calS_{\sph}(X_P)$ is a free right $\calH_{\sph}(M)$-module generated by $f_P$.
\econj
In the case when $P$ is a Borel subgroup this conjecture is proved in \cite{BK}.

\ssec{}{Iwahori part and K-theory}
Let $\calN_{G^{\vee}}$ (resp. $\calN_{M^{\vee}}$) denote the nilpotent cone in the Lie algebra of $G^{\vee}$ (resp. in the Lie algebra of $M^{\vee}$). Let also $\calB_{G^{\vee}},\calB_{M^{\vee}}$ denote the corresponding flag varieties. The cotangent bundle $T^*\calB_{G^{\vee}}$ maps naturally to $\calN_{G^{\vee}}$. Thus we can define
$$
\mathbf{St}_{G^{\vee},M^{\vee}}=T^*\calB_{G^{\vee}}\underset{\calN_{G^{\vee}}}\times T^*\calB_{M^{\vee}}.
$$
This variety is acted on by the group $G^{\vee}\times \CC^{\x}$ (where the second factor acts on $\calN_{G^{\vee}}$ by the formula $t(x)=t^2x$
where $t\in \CC^{\x}$ and $x\in \calN_{G^{\vee}}$. Thus we can consider the complexified equivariant $K$-theory
$K_{M^{\vee}\times \CC^{\x}}(\mathbf{St}_{G^{\vee},M^{\vee}})$. This is a vector space over $\CC[v,v^{-1}=K_{\CC^{\x}}(pt)$. It is easy to see (generalizing the standard construction of ...) that it has a structure of

\end{document}